\documentclass[sn-mathphys]{sn-jnl}
\usepackage[english]{babel}
\jyear{2021}%

\theoremstyle{thmstyleone}%

\newtheorem{theorem}{Theorem}[section]

\newtheorem{example}[theorem]{Example}

\newcommand{\mS}{\mathcal{S}}

\newcommand{\mP}{\mathcal{P}}

\newcommand{\K}{\mathcal{K}}

\newcommand{\bR}{\mathbb{R}}

\newcommand{\on}{\operatorname}

\newcommand{\mpS}{\mathcal{S}_{\on{mp}}}

\begin{document}

\title[Making the computation of {IFS} and {GIFS} deterministic algorithms tractable]{Making the computation of approximations of invariant measures and its attractors for {IFS} and {GIFS}, through the deterministic algorithm, tractable}

\author*[1,2]{\fnm{Rudnei} \sur{D. da Cunha}}\email{rudnei.cunha@ufrgs.br}

\author[1,3]{\fnm{Elismar} \sur{R. Oliveira}}\email{elismar.oliveira@ufrgs.br}

\affil*[1]{\orgdiv{Department of Pure and Applied Mathematics}, \orgname{Instituto de Matem\'atica e Estat\'istica, Universidade Federal do Rio Grande do Sul}, \orgaddress{\street{Av. Bento Gon\c calves 9500}, \city{Porto Alegre}, \postcode{91500-900}, \state{RS}, \country{BRAZIL}}}
\affil[2]{ORCID: 0000-0003-3057-7882}
\affil[3]{ORCID: 0000-0003-2611-0489}

\abstract{
We present algorithms to compute approximations of invariant measures and its attractors for IFS and GIFS, using the deterministic algorithm in a tractable way, with code optimization strategies and use of data structures and search algorithms. The results show that these algorithms allow the use of these (G)IFS in a reasonable running time.
}



\maketitle

\section{Introduction}
\label{sec:intro}
The Hutchinson--Barnsley theory for IFS seek to establish the existence of an unique attractor set, invariant by the fractal operator, and an unique measure of  probability with support on the attractor, invariant by the Markov operator. Those objects are of paramount importance in the extensive applications in several fields of pure and applied sciences.

There are mainly three types of algorithms to approximate the attractor and the invariant measure. The deterministic one, where and initial set or measure is iterated by the respective operators approximating the attractor or the invariant measure w.r.t. the appropriated topology, see \cite{Bar, BDEG} for classical IFS and \cite{JMS16} for GIFS. The discrete one, is similar to the deterministic but an initial step is to introduce a discrete version of the space, an $\varepsilon$-net, and a discrete version of the operator, which is then iterated in the same fashion as the deterministic one producing a discrete set close to the attractor and a discrete measure which is close to the invariant one, \cite{GN, GMN, COS21}. For the last, we have several variations of the original chaos game algorithm introduce by M. Barnsley, \cite{Bar}, where an initial point is iterated choosing, according to some probability, the function to be used. As the process goes on its orbit approximates the attractor set.

As can be seen in \cite{MM3, GN, GMN} and \cite{COS20,COS21} the discrete algorithms exhibit a good performance, but require a lot of technical detail to its implementation. On the other hand deterministic algorithms are easy to describe and implement, in its naive version, but they are impractical computationally. Our aim is to improve this feature.

As described in the sequel, we are interested in deterministic algorithms for IFS and GIFS that require a search for a given value over very large sets of points. A case where it occurs is the approximation of the invariant measure for an IFS or GIFS by iterating the Markov operator from a single point of mass. This task is particularly hard for Idempotent IFS and GIF due its nature.

There are many search algorithms that are tailored to perform geometric scans of a region to locate an specific point, and they all use some form of hierarchical partitioning of the region. In our case, algorithms that deal with discrete representations of the points (for instance, line or column ordering of pixels) are not suited. Instead, we have opted to use a rectangular, regular partition of the region into four quadrants, which leads to the use of quadtrees, a well-known data structure (see \cite{Finkel1974}, \cite{HS1984}, \cite{HS1988} and \cite{SHP2011}, to cite just a few). This is coupled with a fast indexing function to locate quadrants of ever diminishing area, whose limits bracket the point being sought for.

The paper is organized as follows. Section \ref{sec:HB_theory} recalls the Hutchinson--Barnsley theory. In Section \ref{sec:ifs-d} we state the deterministic algorithm for IFS, followed by the introduction of our quadtree-based search algorithm on Section \ref{sec:quadtree}. The natural extension of these algorithms to GIFS is given on Section \ref{sec:gifs-d}. Finally, we conclude with our remarks on Section \ref{sec:conclusion}.

\section{The Hutchinson--Barnsley theory}
\label{sec:HB_theory}

For the convenience of the reader, we will now recall a few basic facts on Iterated Function Systems (IFS for short) and Generalized Iterated Function Systems (GIFS for short).\\ Let $(X,d)$ be a complete, Haussdorf  metric space. By an IFS with probabilities  we mean a triple $\mS=(X,(\phi_j)_{j=1}^L,(p_j)_{j=1}^L)$ so that $\phi_j: X \to X$ and $(X,(\phi_j)_{j=1}^L)$ is an IFS and $p_1,...,p_L \geq 0$ with $\sum_{j=1}^Lp_j=1$. Since IFS are widely known in the literature we will avoid to repeat those definitions for GIFS because they are almost equal, except by the fact that $\phi_j: X\times X \to X$.  Each IFS $\mS=(X,(\phi_j)_{j=1}^L)$ generates the Hutchinson--Barnsley operator $F_\mS:\K(X)\to\K(X)$, where $\K(X)$ is the set of nonempty compact sets of $X$, defined by
$ \displaystyle \forall_{K\in\K(X)}\;F_\mS(K):=\bigcup_{j=1}^L\phi_j(K).$
A set $A_\mS\in\K(X)$ is called the attractor of the IFS $\mS$, if
$A_\mS=F_\mS(A_\mS)$ and for every $K\in\K(X)$, the sequence of iterations $F_\mS^{(n)}(K)\to A_\mS$ w.r.t. the Hausdorff metric.\\
Each IFSp generates also the map $M_\mS:\mP(X)\to\mP(X)$, called the Markov operator, which adjust to every $\mu\in\mP(X)$, the measure $M_\mS(\mu)$ defined by $ \displaystyle M_{\mathcal{S}}(\mu)(B)=\sum_{j=1}^{L} p_{j} \mu\left(\phi_{j}^{-1}(B)\right)$, for any Borel set $B\subset X$.\\
By an invariant measure of an IFSp $\mS$ we mean a (necessarily unique) measure $\mu_\mS\in\mP(X)$ which satisfies $\mu_\mS=M_\mS(\mu_\mS)$ and such that for every $\mu\in\mP(X)$, the sequence of iterates $M^k_\mS(\mu)$ converges to $\mu_\mS$ with respect to the Monge--Kantorovich distance.
The Markov operator $M_\mS$, is also characterized by
\begin{equation}
\int_{X}f \;dM_\mS(\mu)=\sum_{j=1}^Lp_j\int_Xf\circ \phi_j\;d\mu,
\end{equation}
for every IFSp $\mS$ and every continuous map $f:X\to\bR$.
The following result is known (see, for example \cite[Section 4.4]{Hut}).
\begin{theorem}
Each   IFSp on a complete metric space consisting of Banach contractions admits an invariant measure.
\end{theorem}
The Lemma 5.1, from \cite{COS21}, is the basis for the deterministic algorithm, Algorithm~\ref{alg:ifs-d}. Starting with an initial measure of probability $\mu=\delta_{x_0} \in \mP(X)$, each iteration produces a new measure $M_{\mS}(\mu), M_{\mS}^{2}(\mu), ...$, converging to the invariant measure, whose weight in each point of the support, described below, requires to compute the weights of all points which has the same image. Additionally, the set of supporting points describe the deterministic algorithm to approximate de attractor. If for some $N\geq 1$ we have $M_{\mathcal{S}}^{N}(\mu)=\sum_{i=1}^m v_i  \delta_{y_i}\in \mP(X)$, that is, each $v_i\geq 0$ and $\sum_{i=1}^m v_i=1$,  then \begin{equation}\on{supp}(M_{\mathcal{S}}^{N+1}(\mu))=
\{\phi_j(y_{i}):j=1,...,L,\; i\in\{1,...,m\}\}\end{equation} and enumerating this set by $\{z_1,...,z_{m'}\}$, we have:
\begin{equation}\label{update_rule_IFS}
M_{\mathcal{S}}^{N+1}(\mu)=\sum_{r=1}^{m'} v'_r  \delta_{z_r},
\end{equation}
where $\displaystyle v'_{r}=\sum_{\phi_j(y_{i})=z_r} p_j v_{i}$.\\
In other words, the discrete algorithm starts with a product set $D_{0}=\{(x_0,1)\}$ which is the initial value and, at each iteration, the set  $D_{N}=\{(y_1, v_1),...,(y_m, v_m)\}$ is updated by the application of the Markov operator producing the new set $$D_{N+1}=\{(z_1, v'_1),...,(z_{m'}, v'_{m'})\},$$ obtained by the updating rule \eqref{update_rule_IFS}. The first coordinates of $D_{N}$, given by $\on{supp}(M_{\mathcal{S}}^{N+1}(\mu))$ approximate the attractor set $A_{\mathcal{S}}$, and the second coordinates $\{v'_1,...,v'_{m'}\}$, gives the value at each point of the discrete probability $M_{\mathcal{S}}^{N}(\mu)$ approximating the invariant probability $\mu_{\mathcal{S}}$.\\
For a GIFS the updating rule is almost the same, if $M_{\mathcal{S}}^{N}(\mu)=\sum_{i=1}^m v_i  \delta_{y_i}\in \mP(X)$, that is, each $v_i\geq 0$ and $\sum_{i=1}^m v_i=1$,  then \begin{equation}\on{supp}(M_{\mathcal{S}}^{N+1}(\mu))=
\{\phi_j(y_{i_0}, y_{i_1}):j=1,...,L,\;i_0, i_{1} \in\{1,...,m\}\}\end{equation} and enumerating this set by $\{z_1,...,z_{m'}\}$, we have:
\begin{equation}\label{update_rule_GFS}
M_{\mathcal{S}}^{N+1}(\mu)=\sum_{r=1}^{m'} v'_r  \delta_{z_r},
\end{equation}
where $\displaystyle v'_{r}=\sum_{\phi_j(y_{i_0}, y_{i_1})=z_r} p_j v_{i_0}v_{i_1}$.\\
In \cite{MZ} and more recently in \cite{dacunha2021fuzzyset} there was considered the following version in the context of idempotent measures.
Let {$\mathbb{R}_{{\rm max}}:=\mathbb{R} \cup \{-\infty\}$} be the extended set of real numbers. Consider the operations $x \oplus y = \max \{ x, y\}$ and $x \odot y = x+y$. Then we define the \emph{max-plus semiring} $S$ as the algebraic  structure $S=(\mathbb{R}_{{\rm max}}, \oplus, \odot)$. A functional $\mu:C(X) \to \mathbb{R}$ satisfying
  \begin{enumerate}
    \item $\mu(\lambda)=\lambda$ for all $\lambda\in\bR$ (normalization);
    \item $\mu(\lambda \odot \psi)=\lambda \odot \mu(\psi)$, for all $\lambda \in\mathbb{R}$ and $\psi \in C(X)$;
    \item $\mu(\varphi \oplus \psi)=\mu(\varphi) \oplus \mu(\psi)$, for all $\varphi, \psi \in C(X)$,
  \end{enumerate}
  is called an idempotent probability measure (or Maslov measure), \cite{Z, Zai}. A key idea is the \emph{density} of an idempotent probability measure introduced in \cite{Kol88}.
If $\lambda:X\to[-\infty,0]$ is upper semicontinuous and $\lambda(x)=0$ for some $x\in X$, then the map
$ \displaystyle\mu_{\lambda}=\bigoplus_{x\in X}\lambda(x)\odot\delta_x$
is an idempotent measure, that is, $\mu_{\lambda}\in I(X)$. The density $\lambda_{\mu}$ of $\mu\in I(X)$ is uniquely determined.\\
\emph{{Let $\mS=(X,(\phi_j)_{j=1}^m)$ be an IFS and $(q_j)_{j=1}^L$ is a family of real numbers so that $ q_j \leq 0$ for $j=1,...,L$  and,  $\displaystyle\bigoplus_{j=1,...,L} q_j =0$. Then we call the triple $\mpS=(X,(\phi)_{j=1}^L,(q_j)_{j=1}^L)$ as}} a max-plus normalized IFS (which is the idempotent analogous of the IFSp).
Each max-plus normalized IFS $\mpS=(X,(\phi_j)_{j=1}^L,(q_j)_{j=1}^L)$ generates the map $M_{\mathcal{S}}:I(X)\to I(X)$, called as the idempotent Markov operator, which adjust to every $\mu\in I(X)$, the idempotent measure $M_{\mathcal{S}}(\mu)$ defined by:
\begin{equation*}
M_{\mathcal{S}}(\mu):=\bigoplus_{j=1}^Lq_j\odot (I(\phi_j)(\mu))
\end{equation*}
that is, for every $\psi\in C(X)$,
$M_{\mathcal{S}}(\mu)({\psi})=\bigoplus_{j=1}^{L} q_{j}  \odot \mu({\psi}\circ\phi_j)$, because $I(\phi_j)(\mu)(\psi):= \mu({\psi}\circ\phi_j)$.\\
By an invariant idempotent measure \emph{of a {max-plus normalized IFS $\mpS$} we mean  the unique measure $\mu_{\mathcal{S}}\in I(X)$ which satisfies
\begin{equation*}
\mu_{\mathcal{S}}=M_{\mathcal{S}}(\mu_{\mathcal{S}})
\end{equation*}
and such that for every $\mu\in I(X)$, the sequence of iterates $M^{(n)}_{\mathcal{S}}(\mu)$ converges to $\mu_{\mathcal{S}}$ with respect to the $\tau_p$ topology on $I(X)$.\\
We say that a max-plus{ normalized IFS $\mpS$ is}} Banach contractive, if the underlying IFS ${\mS}$ is contractive. The main result of \cite{MZ}, that is \cite[Theorem 1]{MZ}, states that:
\begin{theorem}\label{thm:existence_idemp_measure}
Each Banach contractive {max-plus normalized IFS $\mpS$} on a complete metric space generates the unique invariant idempotent measure $\mu_{\mathcal{S}}$.
\end{theorem}
We notice that Theorem~\ref{thm:existence_idemp_measure} was extended in several ways in our recent works \cite{dacunha2021fuzzyset} and  \cite{dacunha2021existence}.

Now we give a description of the Idempotent Markov operator acting on a finite measure, which is an application of Lemma 5.5 from \cite{dacunha2021fuzzyset}. For $\mu=\bigoplus_{x\in X}\lambda(x)\odot\delta_x\in I(X)$, we have that
\begin{equation*}
  M_{\mathcal{S}}(\mu)=\bigoplus_{y\in X}\lambda_\mS(y)\odot\delta_y
\end{equation*}
where
\[\lambda_\mS(y)=\left\{\begin{array}{cc}\max\{q_j+\lambda(x):j=1,...,L,\;x\in\phi^{-1}_j(y)\},&\mbox{if}\;y\in\bigcup_{j=1}^L\phi_j(X)\\
-\infty,&\mbox{otherwise}\end{array}\right..
\]

Analogously to the classic IFS, if we start with  a singleton $\mu= \delta_{x_{0}} \in I(X)$, the deterministic method consists in to iterate $M_{\mathcal{S}}(\mu), M_{\mathcal{S}}^{2}(\mu),... $ which are also discrete measures with support in a discrete set, converging to the invariant one: \\
In the above frame, if for some $N\geq 1$ we have $M_{\mathcal{S}}^{N}(\mu)=\bigoplus_{i=1}^m v_i \odot\delta_{y_i}\in I(X)$, that is, each $v_i\leq 0$ and $\bigoplus_{i=1}^m v_i=0$,  then \begin{equation}\label{supp_idemp_IFS}\on{supp}(M_{\mathcal{S}}^{N+1}(\mu))=
\{\phi_j(y_{i}):j=1,...,L,\;i=\{1,...,m\}\}\end{equation} and enumerating this set by $\{z_1,...,z_{m'}\}$, we have:
\begin{equation}\label{update_rule_Idemp_IFS}
M_{\mathcal{S}}^{N+1}(\mu)=\bigoplus_{r=1}^{m'} v'_r \odot \delta_{z_r},
\end{equation}
where $\displaystyle v'_{r}=\max_{\phi_j(y_{i})=z_r} q_j + v_{i}$.\\
In other words, the discrete algorithm starts with a direct product set $D_{0}=\{(x_0,1)\}$ which is the initial value and, at each iteration, the set  $D_{N}=\{(y_1, v_1),...,(y_m, v_m)\}$ is updated by the application of the idempotent Markov operator producing the new set $$D_{N+1}=\{(z_1, v'_1),...,(z_{m'}, v'_{m'})\},$$ obtained by the updating rule \eqref{update_rule_Idemp_IFS}. The first coordinates of $D_{N}$, given by $\on{supp}(M_{\mathcal{S}}^{N+1}(\mu))$ approximate the attractor set $A_{\mathcal{S}}$, and the second coordinates $\{v'_1,...,v'_{m'}\}$, gives the value at each point of the discrete probability $M_{\mathcal{S}}^{N}(\mu)$ approximating the invariant idempotent probability $\mu_{\mathcal{S}}$.\\
For an Idempotent GIFS, if $M_{\mathcal{S}}^{N}(\mu)=\bigoplus_{i=1}^m v_i \odot\delta_{y_i}\in I(X)$, that is, each $v_i\leq 0$ and $\bigoplus_{i=1}^m v_i=0$,  then
\begin{equation}\label{supp_idemp_GIFS}\on{supp}(M_{\mathcal{S}}^{N+1}(\mu))=
\{\phi_j(y_{i_0},y_{i_1}):j=1,...,L,\;i_0,i_{1}\in \{1,...,m\}\}\end{equation} and enumerating this set by $\{z_1,...,z_{m'}\}$, we have:
\begin{equation}\label{update_rule_Idemp_GIFS}
M_{\mathcal{S}}^{N+1}(\mu)=\bigoplus_{r=1}^{m'} v'_r \odot \delta_{z_r},
\end{equation}
where $\displaystyle v'_{r}=\max_{\phi_j(y_{i_0}, y_{i_1})=z_r} q_j + v_{i_0} + v_{i_1}$.\\
For additional facts on Hutchinson--Barnsley theory for IFS see \cite{Bar}. For additional facts on idempotent IFS see \cite{Z,MZ,Zai}.  We will not state the analogous facts for GIFS to avoid repetitions; these can be found in \cite{GMM, Mi, MM, MM1}.

\section{Deterministic Algorithm for IFS}
\label{sec:ifs-d}
The standard deterministic algorithm to compute an approximation of both the attractor and the associated invariant measure for an IFS, or an Idempotent IFS, is given in the Algorithm \ref{alg:ifs-d}. To avoid extra technicalities we will work with 2D contained in a rectangle $[a;b]\times [c;d]$.

A number $N$ of iterations is set and a set of $L$ functions $\phi_i: \mathbb{R}^2\rightarrow\mathbb{R}^2$, $1\leq i\leq L$ which describe the attractor is given. The iterations are started over a single point with coordinates $(x,y)$ such that $a\leq x\leq b$ and $c\leq y\leq b$ and some property $p$ associated with that point (say, $p=0$). The output is an array $D$ of triplets $(x,y,p)$. The number of triplets stored in $D$ is referred to by $D.n$.

\begin{algorithm}[H]
\caption{Deterministic Algorithm for IFS}\label{alg:ifs-d}
\begin{algorithmic}[1]
\Function{Deterministic\_IFS}{\textbf{input:} $N,L,x,y,p$; \textbf{output:} $D$}
\Require{$N$, number of iterations}
\Require{$L$, number of $\phi_i$ functions}
\Require{$(x,y)$, coordinates of initial point}
\Require{$p$, some property of the initial point}
\Ensure{$D$, array of $(x,y,p)$ triplets}

\hspace*{-3.6em}\textbf{Local:} $T$, array of $(x,y,p)$ triplets
\State $D.n \leftarrow 1$
\State $D[D.n]\leftarrow (x,y,p)$
\For{$i\leftarrow 1$ \textbf{to} $N$}
    \State $T.n\leftarrow 0$ // initialize the number of triplets stored in $T$
    \For{$j\leftarrow 1$ \textbf{to} $L$}
        \For{$k\leftarrow 1$ \textbf{to} $D.n$}
            \State $(u,v)\leftarrow \phi_{j}(D[k].x,D[k].y)$
            \State \textbf{initialize} $r$ with some appropriate value
            \State \textbf{search} for $(u,v)$ in $T$\label{alg:ifs-search-point}
            \If{$(u,v)$ is not in $T$}
                 \State $T.n\leftarrow T.n+1$
                 \State $T[T.n]\leftarrow (u,v,r)$
            \Else
                \State // $m$ is the index in $T$ such that $T[m].(x,y)=(u,v)$
                \State \textbf{update} the value of $p$ in the triplet $T[m]$ (using $r$)\label{alg:ifs-p-update}
            \EndIf
        \EndFor
    \EndFor
    \State $D \leftarrow T$
\EndFor
\EndFunction
\end{algorithmic}
\end{algorithm}
The algorithm is written as it would be suitable to compute the invariant measure whose support is the attractor. To this end, a search (line \ref{alg:ifs-search-point} of Algorithm \ref{alg:ifs-d}) for a point $(u,v)$ that has been produced from a point $D[k]$ produced in the previous iteration, with an associated property $r$, is made over the points being produced in the current iteration (which are stored on a local array $T$ of triplets $(x,y,p)$ similar to $D$). If $(u,v)$ is found on $T$, then its associated property $p$ is modified, dependent on how the measure is defined: that could be expressed as the maximum between, or a sum involving $p$ and $r$.

We note that $p$ is updated (Algorithm \ref{alg:ifs-d}, line \ref{alg:ifs-p-update}) according to Equation (\ref{update_rule_IFS}) (for classic IFS) and to Equation (\ref{update_rule_Idemp_IFS}) (for idempotent IFS). For the latter, the variable $r$ is initialized with $p_{j}+D[k].p$ and the update is $T[m].p\leftarrow \max(T[m].p,r)$.

 Of course, if one does not wish to compute such a measure, then the algorithm becomes much simpler: there is no need to perform any searches and it is just a matter of applying the $\phi_i$ functions over the set of points produced at the previous iteration.

 The main issue arising with Algorithm \ref{alg:ifs-d} (or, indeed, with its simpler version) is that the number of points produced at each iteration grow considerably, to a point where it becomes impracticable to execute it in a reasonable time, or simply the amount of memory required to store the points in the arrays $D$ and $T$ is not available. Another issue, one that becomes much more noticeable as the iterations proceed, is the number of operations required to search for $(u,v)$ in $T$.

\section{Optimizing Algorithm \ref{alg:ifs-d} using a better search algorithm}
\label{sec:quadtree}
As seen above, the use of a linear search leads to an overall cost in terms of the number of searches that is exponential in nature. It is not possible to reduce this exponential growth, because it is the very essence of Algorithm \ref{alg:ifs-d}; however, we can make it more tractable if a better search algorithm is employed. This may allow, for instance, to perform one more iteration in Algorithm \ref{alg:ifs-d}, and that may be just enough to obtain a good image of the attractor, or a more refined value for the measure being computed.

Since an IFS produces an image that is fractal in nature and it has a contractivity property, the attractor is a compact set, meaning that the points generated along the iterations will all fall within a suitable rectangular region $X=[a;b]\times [c;d]$ in $\mathbb{R}^2$ such that any point with coordinates $(x,y)$ along the Cartesian axes satisfy $a\leq x\leq b$ and $c\leq y\leq b$, a hierarchical subdivision of the rectangular region is well indicated to help one to obtain a faster search algorithm. This leads to the use of a quadtree to organize the points and ease searching for one, as the search can be done only on the regions where a particular point resides, discarding all the others.

\begin{figure}[H]
\centering
\includegraphics[width=5cm]{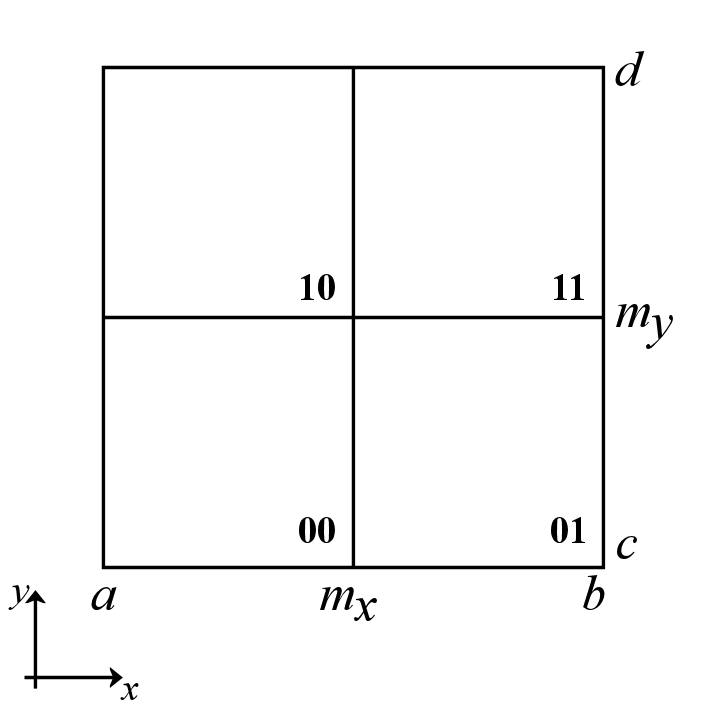}
\caption{Dividing a region into four quadrants; each quadrant is uniquely identified by a two-bit Gray code.}
\label{fig:division-quadrants}
\end{figure}

Given the intervals $[a;b]$ and $[c;d]$, we proceed by computing their midpoints $m_x$ and $m_y$. With these six values the limits of four quadrants can be defined, as in Figure \ref{fig:division-quadrants}. Each quadrant is then assigned a two-bit Gray code, in which the code between two neighbouring quadrants (i.e. those that share a common edge) differs by just one bit. Thus, the quadrants are enumerated from $0$ to $3$.

We now define the function $Q(x,y)$ which gives the quadrant number as follows:
\begin{eqnarray}
Q(x,y)&=&2i(y,c,d)+i(x,a,b)\label{eq:q}\\
i(v,\alpha,\beta)&=&\left\{\begin{array}{cc}0,&\frac{v-\alpha}{\beta-\alpha}<\frac{1}{2}\\1,&\emph{otherwise}\end{array}\right., \alpha\leq v\leq \beta
\end{eqnarray}
where $0\leq Q(x,y)\leq 3$. Therefore, any point within the region may be attributed to a unique quadrant (we may also refer to the value returned by $Q(x,y)$ in terms of its related two-bit Gray code).

The quadtree data structure is made up of nodes associated to quadrants. Its root node, of course, refers to the region $X$ as defined before. Each node stores the values $a, b, c, d$ along the $x$ and $y$ axes that define the quadrant region. If there are any points within its region, the node will also store an array of integer values $I$ with at most $n_{\max}$ entries. This array contains the indices of the points stored in $T$ (on Algorithm \ref{alg:ifs-d}) which belong to the quadrant.

The node may also hold four node sons, in case it has been divided when trying to insert a point on the quadtree. This will happen whenever a node $K$ was supposed to hold a point $(x,y)$ for which $Q(x,y)=k$, but the array $I_K$ has no more free entries. In this case, the node is subdivided into four sons, by the midpoints of $[a;b]$ and $[c;d]$, as in Figure \ref{fig:division-quadrants}. The points associated to node(quadrant) $K$ are redistributed among its sons and, finally, the point that was being inserted (and caused the subdivision) is assigned to one of its sons, recursively. Only leaf nodes (i.e. nodes without sons) have the array $I$ and any search for a point $(x,y)$ occurs only on a leaf node.

These ideas are presented in Algorithm \ref{alg:quadtree-search-insert}. We also present Algorithm \ref{alg:ifs-d-quadtree}, which is a modification of Algorithm \ref{alg:ifs-d} to use the quadtree data structure.

\begin{algorithm}
\caption{Quadtree search and insert}\label{alg:quadtree-search-insert}
\begin{algorithmic}[1]
\Function{Quadtree\_Search\_and\_Insert}{\textbf{input:} $Q,u,v,r$; \textbf{output:} $D$}
\Require{$Q$, a quadtree node}
\Require{$(u,v)$, coordinates of point to search}
\Require{$r$, some property of the point to search}
\Ensure{$D$, array of $(x,y,p)$ triplets}
\If{node $Q$ has sons}
    \State $i\leftarrow Q(u,v)$
    \State \Call{Quadtree\_Search\_and\_Insert}{$Q.sons[i],u,v,r,D$}
\Else
    \State $f\leftarrow 0$
    \For{$j\leftarrow 1$ \textbf{to} $Q.I.n$}\label{alg:linear-search-quadtree-search-insert}
        \State // linear search for $(u,v)$ on $Q.I$
        \If{$u=D[Q.I[j]].x$ AND $v=D[Q.I[j]].y$}
            \State \textbf{update} the value of $p$ in the triplet $D[Q.I[j]]$ (using $r$)\label{alg:update-p-quadtree}
            \State $f\leftarrow 1$
            \State \textbf{break}
        \EndIf
    \EndFor
    \If{$f=0$}
        \State // if $(u,v)$ was not found on $Q.I$
        \If{$Q.I.n<n_{\max}$}
            \State // if there are free entries available on $Q.I$,
            \State // add $(u,v)$ to $D$ and its index to $Q.I$
            \State $D.n\leftarrow D.n+1$
            \State $D[D.n]\leftarrow (u,v,r)$
            \State $Q.I.n\leftarrow Q.I.n+1$
            \State $Q.I[Q.I.n]\leftarrow D.n$
        \Else
            \State \textbf{divide} node $q$ into four sons
            \For{$j\leftarrow 1$ \textbf{to} $Q.I.n$}
                \State // distribute the indices of $Q$ among its sons
                \State $i\leftarrow Q(D[q.I[j]].x,D[q.I[j]].y)$
                \State $Q.sons[i].I.n\leftarrow Q.sons[i].I.n+1$
                \State $Q.sons[i].I[Q.sons[i].I.n]\leftarrow Q.I[j]$
            \EndFor
            \State // insert $(u,v,r)$ into the appropriate son of $Q$
            \State $i\leftarrow Q(u,v)$
            \State \Call{Quadtree\_Search\_and\_Insert}{$Q.sons[i],u,v,r,D$}
        \EndIf
    \EndIf
\EndIf
\EndFunction
\end{algorithmic}
\end{algorithm}

\begin{algorithm}
\caption{Deterministic Algorithm for IFS with quadtree-based search}\label{alg:ifs-d-quadtree}
\begin{algorithmic}[1]
\Function{Deterministic\_IFS\_quadtree}{\textbf{input:} $N,L,x,y,p$; \textbf{output:} $D$}
\Require{$N$, number of iterations}
\Require{$L$, number of $\phi_i$ functions}
\Require{$(x,y)$, coordinates of initial point}
\Require{$p$, some property of the initial point}
\Ensure{$D$, array of $(x,y,p)$ triplets}

\hspace*{-3.6em}\textbf{Local:} $T$, array of $(x,y,p)$ triplets

\hspace*{-3.6em}\textbf{Local:} $R$, root node of quadtree
\State $D.n \leftarrow 1$
\State $D[D.n]\leftarrow (x,y,p)$
\For{$i\leftarrow 1$ \textbf{to} $N$}\label{alg:iteration-ifs-d-quadtree}
    \State $T.n\leftarrow 0$ // initialize the number of triplets stored in $T$
    \State \textbf{create} root node of quadtree, $R$
    \For{$j\leftarrow 1$ \textbf{to} $L$}
        \For{$k\leftarrow 1$ \textbf{to} $D.n$}
            \State $(u,v)\leftarrow \phi_{j}(D[k].x,D[k].y)$
            \State \textbf{initialize} $r$ with some appropriate value
            \State \Call{Quadtree\_Search\_and\_Insert}{$R,u,v,r,T$}
        \EndFor
    \EndFor
    \State $D \leftarrow T$
\EndFor
\EndFunction
\end{algorithmic}
\end{algorithm}

\begin{figure}[H]
\centering
\includegraphics[width=7cm]{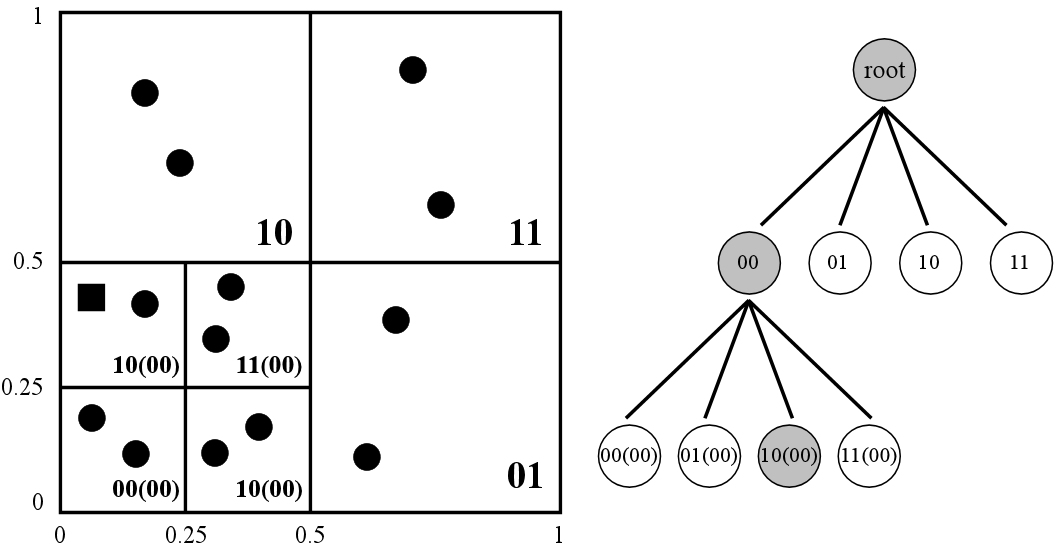}
\caption{Looking for point $(x,y)=(0.1,0.425)$ (black square) on the quadtree; nodes traversed on the quadtree during the search are marked in grey. Quadrant number is indicated in two-bit Gray code; number in brackets indicate the number of the parent quadrant.}
\label{fig:example-quadtree-search-insert}
\end{figure}

We give now an example of the workings of Algorithm \ref{alg:quadtree-search-insert}. Consider a fractal within the region $[0,1]^2$, that at an iteration $i$ the number of points $n_{i}$ generated by Algorithm \ref{alg:ifs-d} (marked as circles) is $n_{i}=13$ and that $n_{\max}=2$, as shown in Figure \ref{fig:example-quadtree-search-insert}. Now suppose the point $(x,y)=(0.1,0.425)$ (marked as a square) is to be looked for on the quadtree shown in the picture, which has height $h=2$.

Since the root node has sons, computing $Q(0.1,0.425)$ gives quadrant $00$ (i.e. son $0$ of the root) as the next node to be traversed on the quadtree. Upon visiting this node, since it also has sons, again $Q(0.1,0.425)$ is computed but this time (since the values of $b$ and $c$ of son $0$ of the root are different from those of the root node) it gives quadrant $10$ (i.e. son $2$ of son $0$ of the root). Now, since this is a leaf node, point $(0.1,0.425)$ is searched for on its array $I$ and either will be found or will be added to $I$ otherwise. Therefore, searching for a point on the quadtree is equivalent to traversing a list of $h$ nodes and then performing a linear search when at most $n_{\max}$ comparisons will be made.

\subsection{Algorithm \ref{alg:ifs-d-quadtree} in practice}
\label{sec:ifs-d-quadtree-in-practice}
To illustrate the mechanics of Algorithm \ref{alg:ifs-d-quadtree}, we consider the classic geometric fractal, Maple Leaf, defined by $L=4$ functions $\phi_i$:
\begin{example}\label{ex:maple-leaf}Maple Leaf:
\[
\left\{
  \begin{array}{ll}
      \phi_1(x,y)=&(0.8x+0.1,0.8y+0.04)\\
      \phi_2(x,y)=&(0.5x+0.25,0.5y+0.4)\\
      \phi_3(x,y)=&(0.355x-0.355y+0.266,0.35x+0.355y+0.078)\\
      \phi_4(x,y)=&(0.355x+0.355y+0.378,-0.355x+0.355y+0.434)
  \end{array}
\right.
\]
on the region $X=[0,1]^2$, with $p_{1}=0$, $p_{2}=-7$, $p_{3}=-3$ and $p_{4}=-7$.
\end{example}

\begin{figure}
  \centering
  \begin{tabular}{c c}
  \includegraphics[width=3.85cm]{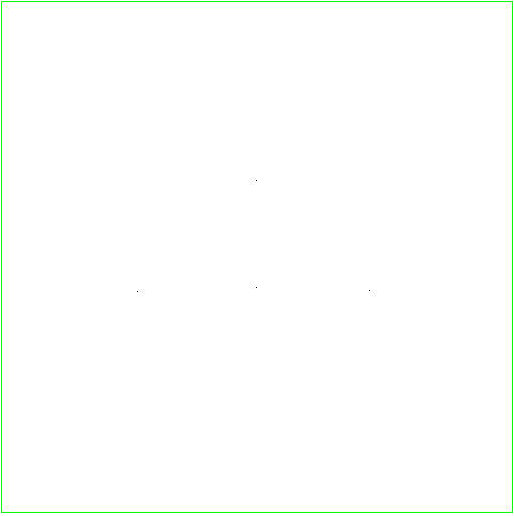}&\includegraphics[width=3.85cm]{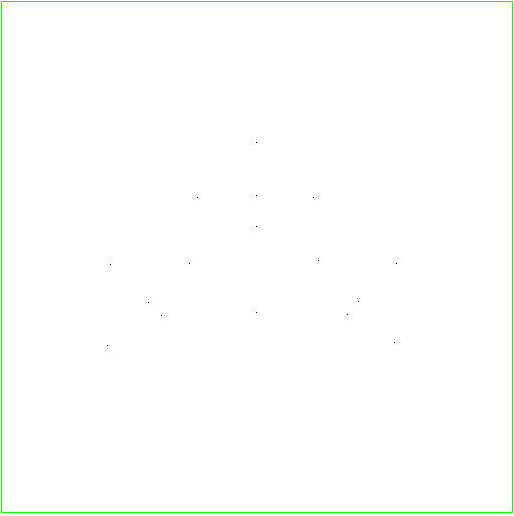}\\
  \includegraphics[width=3.85cm]{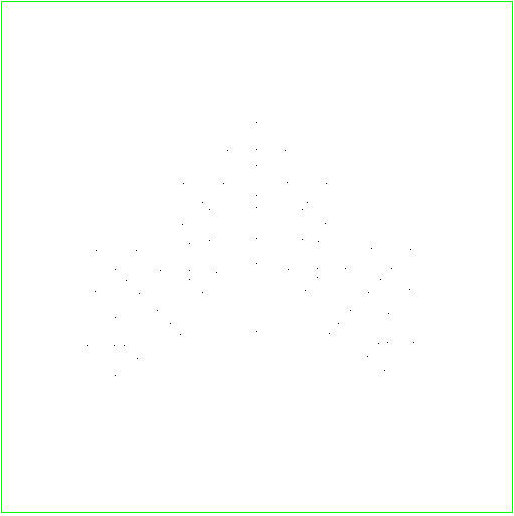}&\includegraphics[width=3.85cm]{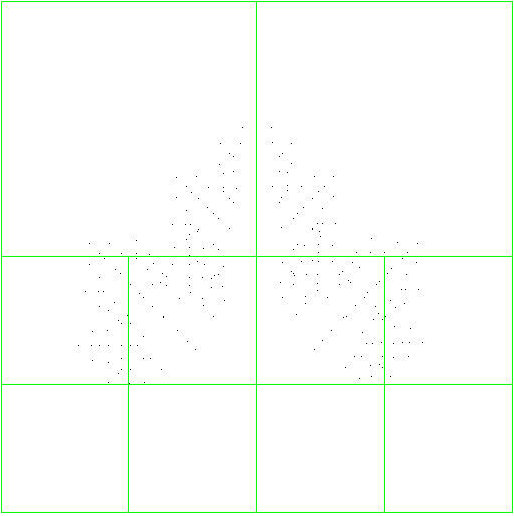}\\
  \includegraphics[width=3.85cm]{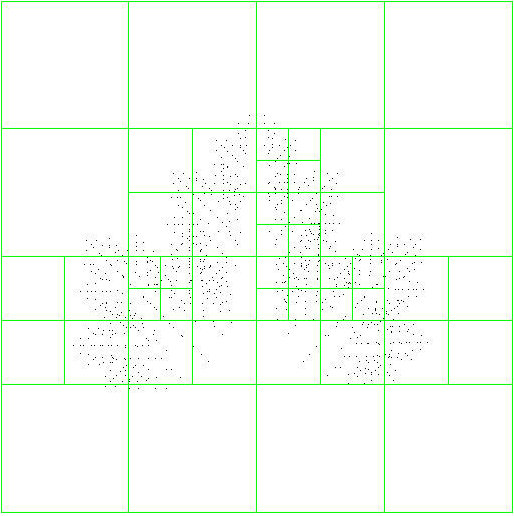}&\includegraphics[width=3.85cm]{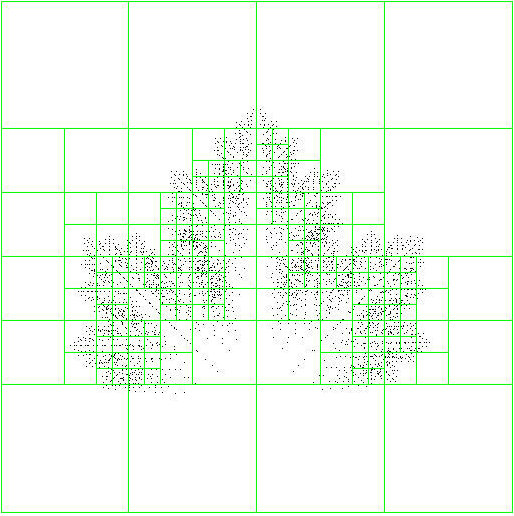}\\
  \includegraphics[width=3.85cm]{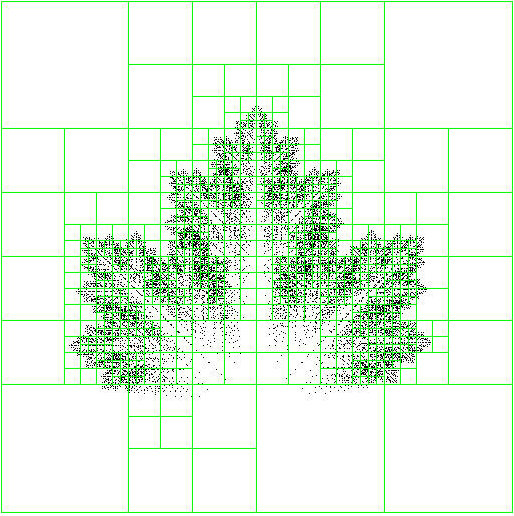}&\includegraphics[width=3.85cm]{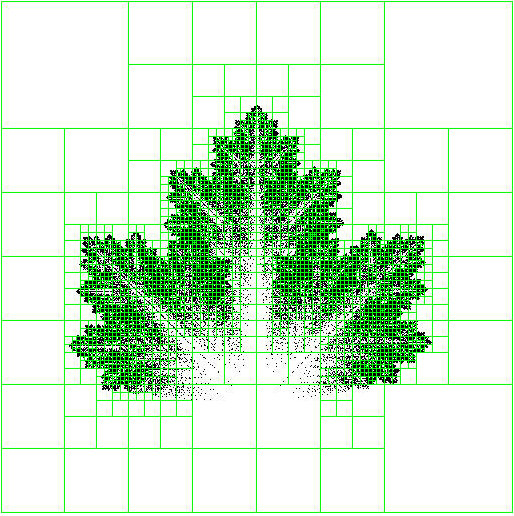}\\
  \includegraphics[width=3.85cm]{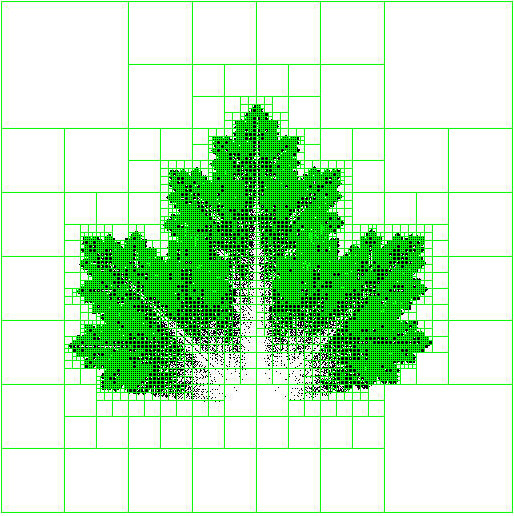}&\includegraphics[width=3.85cm]{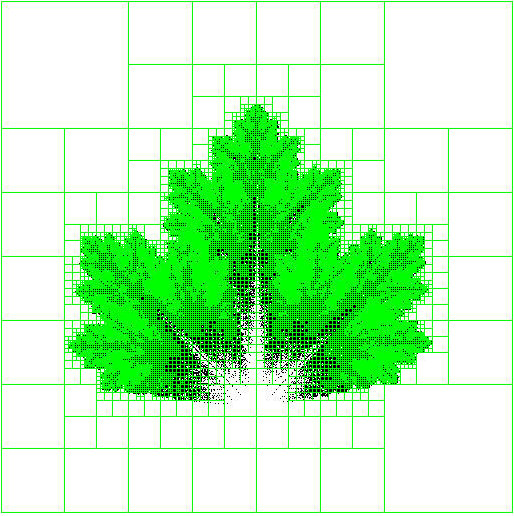}\\
  \end{tabular}
  \caption{The quadtrees placed over the Maple Leaf fractal, $n_{\max}=64$: from top to bottom, left to right, the images ($512\times 512$ pixels) show the points generated from iteration $1$ ($n=4$ points) to $10$ ($n=1048534$ points).}\label{fig:quadtree-maple-leaf}
\end{figure}

Figure \ref{fig:quadtree-maple-leaf} shows the quadtrees placed over the image obtained from the points generated by Algorithm \ref{alg:ifs-d-quadtree}. In this example, we used $n_{\max}=64$. Since the number of points generated at the end of each iteration $i\geq 1$ is $4^i$, only after the fifth iteration there will be node divisions, as can be seen in the images. Note also that there may be points generated in an iteration that are already stored, hence the difference in the number of points after the tenth iteration ($1048534$) instead of the expected $4^{10}=1048576$.

\begin{figure}
  \centering
  \begin{tabular}{c c}
  \includegraphics[width=3.85cm]{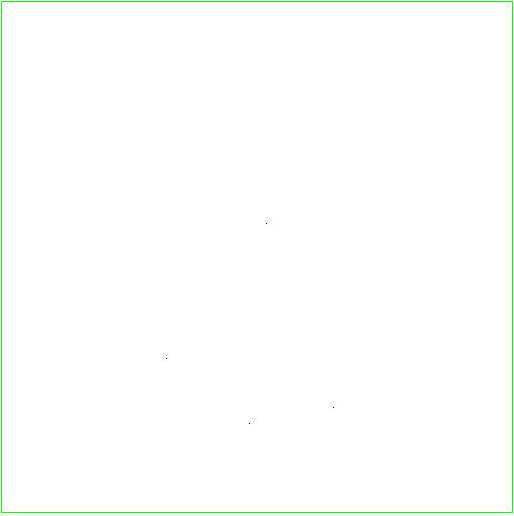}&\includegraphics[width=3.85cm]{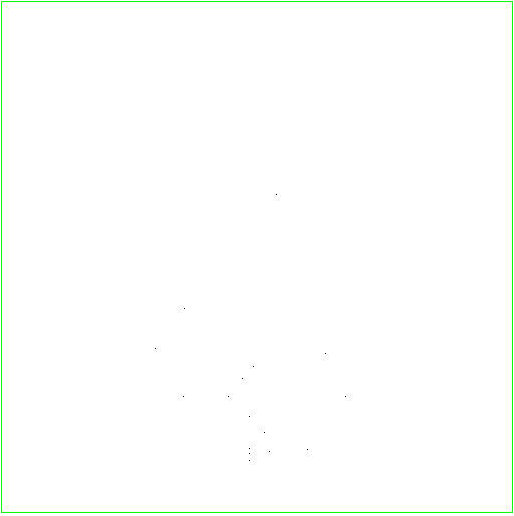}\\
  \includegraphics[width=3.85cm]{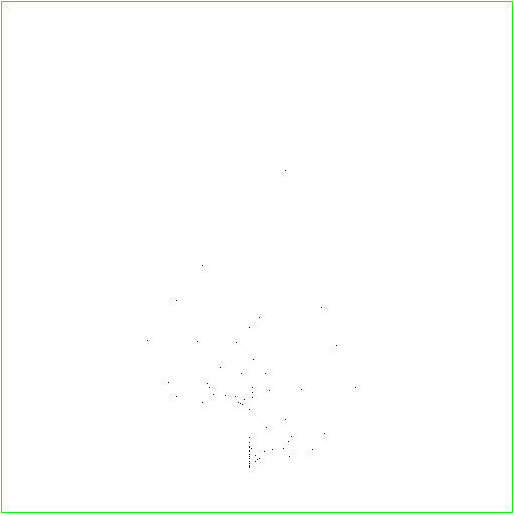}&\includegraphics[width=3.85cm]{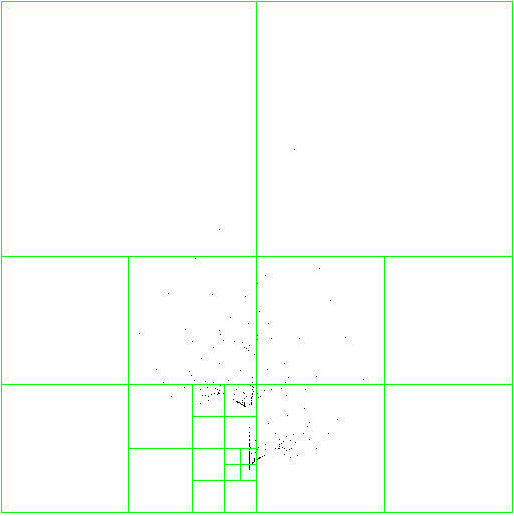}\\
  \includegraphics[width=3.85cm]{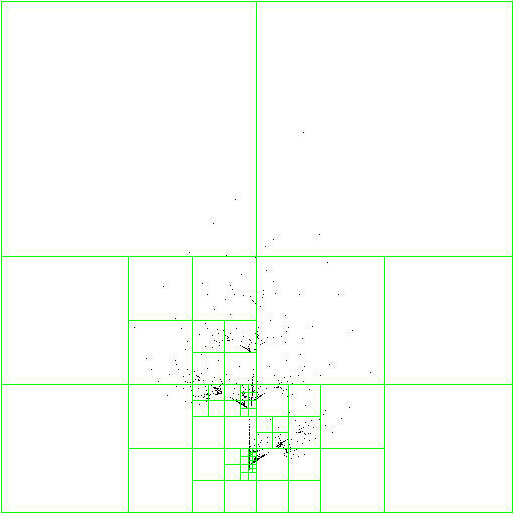}&\includegraphics[width=3.85cm]{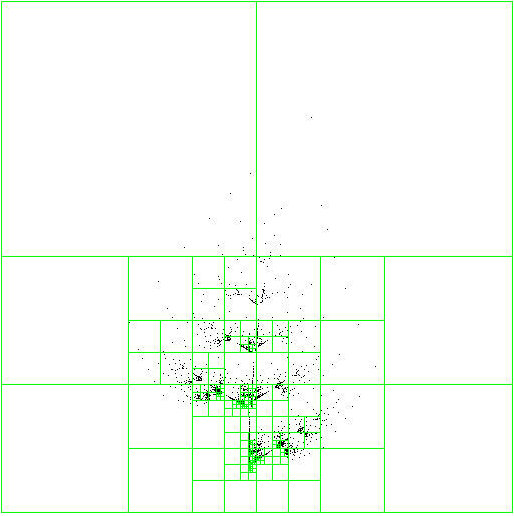}\\
  \includegraphics[width=3.85cm]{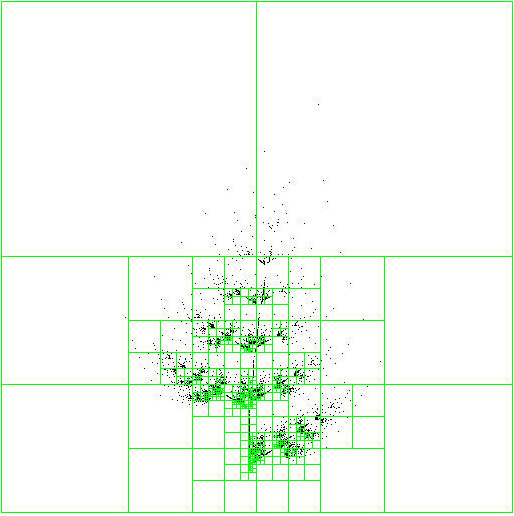}&\includegraphics[width=3.85cm]{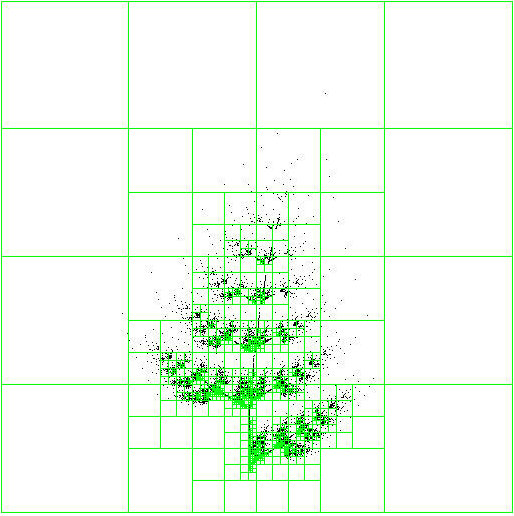}\\
  \includegraphics[width=3.85cm]{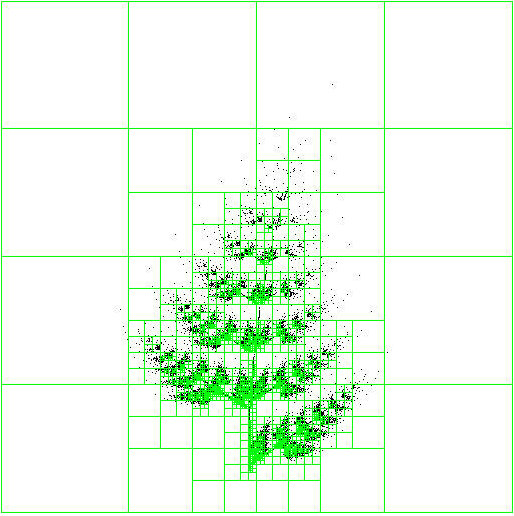}&\includegraphics[width=3.85cm]{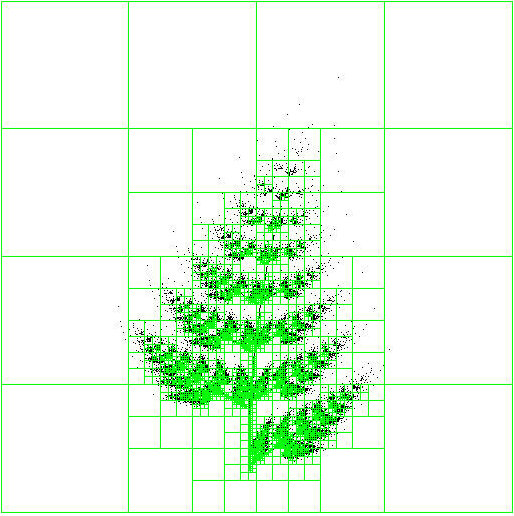}\\
  \end{tabular}
  \caption{The quadtrees placed over the Barnsley Fern fractal, $n_{\max}=64$: from top to bottom, left to right, the images ($512\times 512$ pixels) show the points generated from iteration $1$ ($n=4$ points) to $10$ ($n=1014270$ points).}\label{fig:quadtree-barnsley-fern}
\end{figure}

To llustrate further, we consider an example with a different distribution of points across the region, showing in Figure \ref{fig:quadtree-barnsley-fern} the quadtrees placed over the fractal.
\begin{example}\label{ex:discrete-Fern} %\ero{FILIP: Good!}
This example is based on a very well-known fractal, the Barnsley Fern. It is defined by
\[
\left\{
  \begin{array}{ll}
      \phi_1(x,y)=&(0.856x+0.0414y+0.07,-0.0205x+0.858y+0.147)\\
      \phi_2(x,y)=&(0.244x-0.385y+0.393,0.176x+0.224y+0.102)\\
      \phi_3(x,y)=&(-0.144x+0.39y+0.527,0.181x+0.259y-0.014) \\
      \phi_4(x,y)=&(0.486,0.031x+0.216y+0.05)
  \end{array}
\right.
\]
on the region $X=[0,1]^2$, with $p_{1}=-11$, $p_{2}=-7$, $p_{3}=0$ and $p_{4}=0$.
\end{example}

It should be clear from the previous discussion that $n_{\max}$ plays an important role in Algorithm \ref{alg:quadtree-search-insert}. First, consider that this algorithm traverses the quadtree along a distinct path, visiting only the nodes computed by $Q(u,v)$, until reaching a leaf node (we remember that the height of a (quad)tree is the maximum distance of any node from the root). When this leaf node is visited, then a linear search for $(u,v)$ (of complexity $O\left(n_{\max}\right)$, line \ref{alg:linear-search-quadtree-search-insert} of the algorithm) is made over the points on $D$ that are referred to by the node (indices of entries of $D$ stored on $I$). If $(u,v)$ is found, then $p$ is updated as in Algorithm \ref{alg:ifs-d}.

If the linear search fails, then $(u,v)$ will be stored in $D$ and its index on $D$ is stored on $I$, \emph{if there are available entries on $I$}; otherwise, the node is divided into four sons, the points assigned to it are distributed among its sons and Algorithm \ref{alg:quadtree-search-insert} is called recursively to store $(u,v)$ on one of its sons (which in itself may cause further node divisions).

If $n_{\max}$ is small then, for a given number $n$ of distinct points generated during an iteration (line \ref{alg:iteration-ifs-d-quadtree} on Algorithm \ref{alg:ifs-d-quadtree}), there will be many subdivisions of the nodes on the quadtree, increasing its height, but the linear searches will make few comparison tests. Conversely, a larger $n_{\max}$ value reduces the height of the quadtree, but increases the cost of the linear search on a leaf node.

To ascertain the behaviour of Algorithm \ref{alg:ifs-d-quadtree}, we made a number of runs of our implementation written in \textsc{Fortran 2003} compiled with \textsc{gfortran 10.2.0} with \textsc{-O3} optimization on a computer with an Intel Core i5-6400T 2.20 GHz processor and 6 GB of DDR3 RAM. The examples used were the Maple Leaf defined earlier and the Barnsley Fern.

Table \ref{tab:ifs-d-quadtree} shows the quadtree height, $h$, at the end of $N=10$ iterations, and the execution time (in seconds) of our implementation of Algorithm \ref{alg:ifs-d-quadtree}. For the sake of comparison, we also present the execution time (in seconds) of our implementation of Algorithm \ref{alg:ifs-d}, using a linear search. We note that the number of points generated after $10$ iterations was $n=1048534$ for the Maple Leaf and $n=1014270$ for the Barnsley Fern. It is evident from the data presented that, a) the use of the quadtree provides an execution time that is over $400$ times faster and b) there is an optimal value for $n_{\max}$, namely $64$, for which the least execution time was obtained, among the values used for $n_{\max}$.
\begin{table}
\centering
\caption{Maximum quadtree height and execution times of Algorithm \ref{alg:ifs-d-quadtree}.}\label{tab:ifs-d-quadtree}
\begin{tabular}{c|c|rrr}
  \hline
          & Algorithm \ref{alg:ifs-d} & \multicolumn{3}{c}{Algorithm \ref{alg:ifs-d-quadtree}} \\
     Example & Time [s] & $n_{\max}$ & $h$ & Time [s] \\\hline
  Maple Leaf & 468.600 &    2 & 17 & 1.8610 \\
             &         &    4 & 14 & 1.4590 \\
             &         &    8 & 12 & 1.2660 \\
             &         &   16 & 12 & 1.1920 \\
             &         &   32 & 11 & 1.1120 \\
             &         &   64 & 10 & 1.1000 \\
             &         &  128 &  9 & 1.2200 \\
             &         &  256 &  9 & 1.1730 \\
             &         &  512 &  8 & 1.2900 \\
             &         & 1024 &  8 & 1.5600 \\
  Barnsley Fern & 483.479 &    2 & 26 & 2.3840 \\
                &         &    4 & 25 & 1.7740 \\
                &         &    8 & 24 & 1.4530 \\
                &         &   16 & 23 & 1.2850 \\
                &         &   32 & 22 & 1.2120 \\
                &         &   64 & 21 & 1.1780 \\
                &         &  128 & 20 & 1.1830 \\
                &         &  256 & 19 & 1.2430 \\
                &         &  512 & 17 & 1.4140 \\
                &         & 1024 & 13 & 1.7570 \\
  \hline
\end{tabular}
\end{table}

\begin{figure}
  \centering
  \begin{tabular}{c c}
  \includegraphics[height=3.85cm]{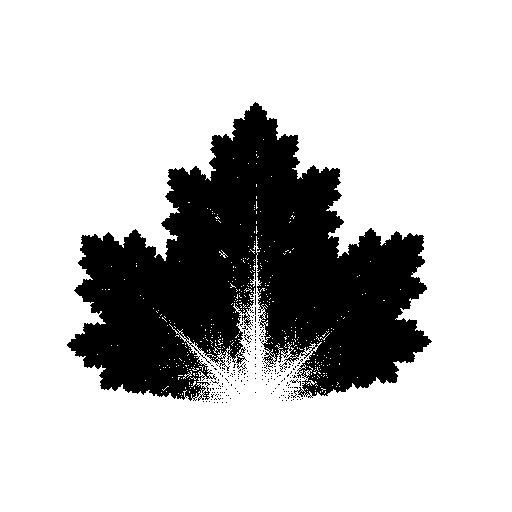}&\includegraphics[height=3.85cm]{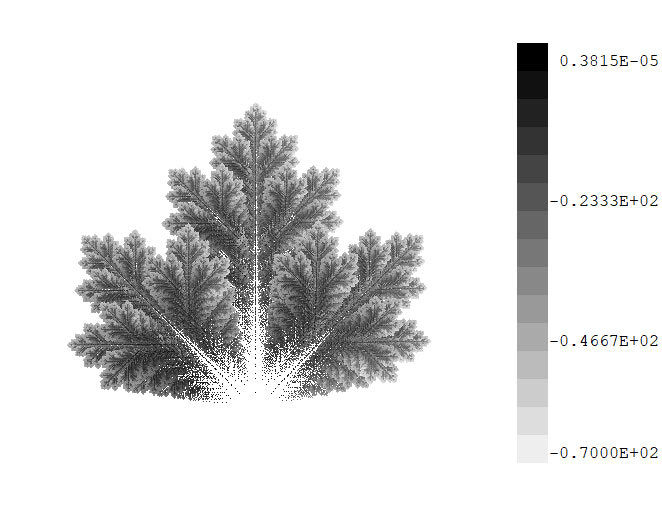}\\
  \includegraphics[height=3.85cm]{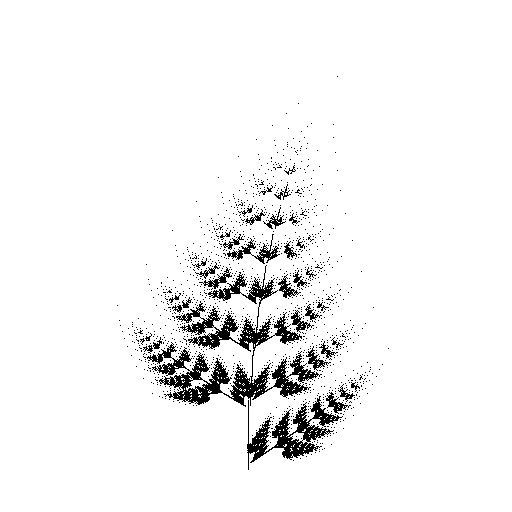}&\includegraphics[height=3.85cm]{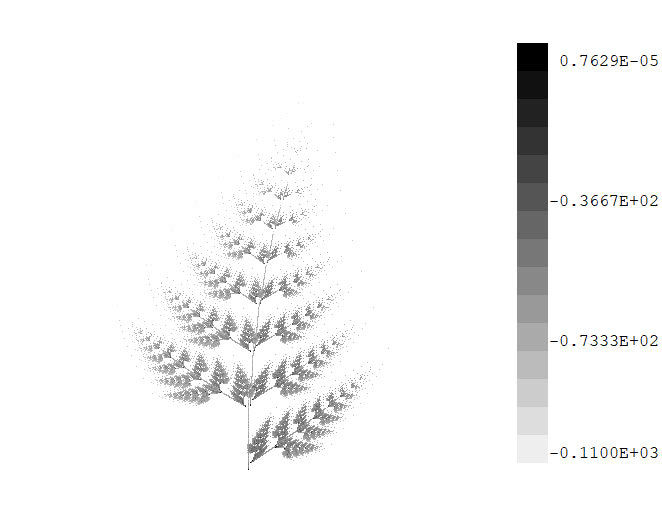}\\
  \end{tabular}
  \caption{The approximation of the attractor and the greyscale representation of the invariant measure for the Maple Leaf (top) and Barnsley Fern (bottom) idempotent IFSs.}\label{fig:ifs-idempotent}
\end{figure}

Finally, to show an example of what we are actually interested in computing with IFSs, we present in Figure \ref{fig:ifs-idempotent} the approximation of the attractor and the greyscale image representing the invariant measure for the Maple Leaf and Barnsley Fern idempotent IFSs. For more details, we refer the reader to \cite[Lemma 5.1]{COS21}.

\section{Complexity analysis of the deterministic IFS algorithms}
\label{sec:complexity-ifs-d}
To establish the complexity of Algorithm \ref{alg:ifs-d}, we introduce the notation $n_{i}$, which is the number of points produced in each iteration, with $n_{0}=1$ (i.e. for the initial point). We assume that the search in line \ref{alg:ifs-search-point} of Algorithm \ref{alg:ifs-d} always fails, meaning that the searched-for point $(u,v)$ will always be stored in $T$ and therefore $n_{i-1}<n_{i}$ (that will give an upper-bound on the number of points generated at each iteration). Typically, one has sequences $n_{i}=L^i$.

Writing the expression for the total number of searches $S_{IFS}$ made in the algorithm, we obtain
\begin{equation}
S_{IFS}=\sum_{i=1}^{N}{\sum_{j=1}^{L}{\sum_{k=1}^{n_{i-1}}{f_{k}}}}
\label{eq:number-searches-ifs-d}
\end{equation}
where $f_{k}$ is the number of comparison tests made on a search over a set of $n_{k}$ values. If a linear search is used, then $f_{k}\in O\left(n_{k}\right)$. In the worst case, $f_{k}=n_{k}$ and Equation (\ref{eq:number-searches-ifs-d}) reduces to
\begin{equation}
S_{IFS}=-{\frac {N{L}^{2}}{L-1}}+\frac{L^{2}}{L-1}\sum_{i=1}^{N}{L^{L^{i-1}}}\simeq-{\frac {N{L}^{2}}{L-1}}+\frac{L^{2}}{L-1}L^{L^{N-1}}
\label{eq:ifs-linear-search}
\end{equation}
where the summation is dropped since the term $L^{L^{N-1}}$ dominates the summation asymptotically as $N\rightarrow \infty$ and, therefore, $S_{IFS}\in O\left(L^{L^{N}}\right)$.

We note that the quadtree-based search has a complexity $O\left(n_{\max}\right)$. This is because the path traversed from the root to a leaf node during a search has at most length $h$ (there may be paths with shorter lengths, it depends on the distribution of points generated along the iterations - see Figure \ref{fig:quadtree-maple-leaf}) and a linear search of at most $n_{\max}$ elements is carried out on a leaf node.

If we choose $n_{\max}=64$ and take $f(k)=n_{\max}$ in Equation (\ref{eq:number-searches-ifs-d}), then Algorithm \ref{alg:ifs-d-quadtree} has a total number of searches given as
\begin{equation}
S_{IFSq}=\frac{64\left({L}^{N+1}-L\right)}{L-1}\in O\left(L^{N}\right).
\label{eq:B-ifs-quadtree}
\end{equation}
and, therefore,
\begin{equation}
S_{IFSq}\ll S_{IFS}.
\label{eq:comparison-ifs-ifsq}
\end{equation}

\section{Extension to Generalized IFS}
\label{sec:gifs-d}
The ideas presented in the previous sections extend naturally to Generalized IFSs. An algorithm to compute the attractor of a deterministic IFS is given in Algorithm \ref{alg:gifs-d}, with a similar notation to that of Algorithm \ref{alg:ifs-d}, and its version using a quadtree is given in Algorithm \ref{alg:gifs-d-quadtree}. They differ from the IFS algorithms in that the $\phi_i$ functions are now $\phi_i: \mathbb{R}^4\rightarrow\mathbb{R}^2$, $1\leq i\leq L$, and also that the set of points at each iteration is obtained by two nested loops of length $D.n$, leading to  $(D.n)^2$ points being produced (at most); this is what makes GIFS costlier to compute than an IFS, since the number of generated points at each iteration grows much more rapidly.

Also, we note that the updates of $p$ on Algorithm \ref{alg:gifs-d}, line \ref{alg:gifs-p-update} and on Algorithm \ref{alg:quadtree-search-insert}, line \ref{alg:update-p-quadtree} are made according to Equation (\ref{update_rule_GFS}) (for classic GIFS) and to Equation (\ref{update_rule_Idemp_GIFS}) (for idempotent GIFS).

\begin{algorithm}
\caption{Deterministic Algorithm for GIFS}\label{alg:gifs-d}
\begin{algorithmic}[1]
\Function{Deterministic\_GIFS}{\textbf{input:} $N,L,x,y,p$; \textbf{output:} $D$}
\Require{$N$, number of iterations}
\Require{$L$, number of $\phi_i$ functions}
\Require{$(x,y)$, coordinates of initial point}
\Require{$p$, some property of the initial point}
\Ensure{$D$, array of $(x,y,p)$ triplets}

\hspace*{-3.6em}\textbf{Local:} $T$, array of $(x,y,p)$ triplets
\State $D.n \leftarrow 1$
\State $D[D.n]\leftarrow (x,y,p)$
\For{$i\leftarrow 1$ \textbf{to} $N$}
    \State $T.n\leftarrow 0$ // initialize the number of triplets stored in $T$
    \For{$j\leftarrow 1$ \textbf{to} $L$}
        \For{$k\leftarrow 1$ \textbf{to} $D.n$}
            \For{$l\leftarrow 1$ \textbf{to} $D.n$}
                \State $(u,v)\leftarrow \phi_{j}(D[k].x,D[k].y,D[l].x,D[l].y)$
                \State \textbf{initialize} $r$ with some appropriate value
                \State \textbf{search} for $(u,v)$ in $T$\label{alg:gifs-search-point}
                \If{$(u,v)$ is not in $T$}
                     \State $T.n\leftarrow T.n+1$
                     \State $T[T.n]\leftarrow (u,v,r)$
                \Else
                    \State // $m$ is the index in $T$ such that $T[m].(x,y)=(u,v)$
                    \State \textbf{update} the value of $p$ in the triplet $T[m]$ (using $r$)\label{alg:gifs-p-update}
                \EndIf
            \EndFor
        \EndFor
    \EndFor
    \State $D \leftarrow T$
\EndFor
\EndFunction
\end{algorithmic}
\end{algorithm}

\begin{algorithm}
\caption{Deterministic Algorithm for GIFS with quadtree-based search}\label{alg:gifs-d-quadtree}
\begin{algorithmic}[1]
\Function{Deterministic\_GIFS\_quadtree}{\textbf{input:} $N,L,x,y,p$; \textbf{output:} $D$}
\Require{$N$, number of iterations}
\Require{$L$, number of $\phi_i$ functions}
\Require{$(x,y)$, coordinates of initial point}
\Require{$p$, some property of the initial point}
\Ensure{$D$, array of $(x,y,p)$ triplets}

\hspace*{-3.6em}\textbf{Local:} $T$, array of $(x,y,p)$ triplets

\hspace*{-3.6em}\textbf{Local:} $R$, root node of quadtree
\State $D.n \leftarrow 1$
\State $D[D.n]\leftarrow (x,y,p)$
\For{$i\leftarrow 1$ \textbf{to} $N$}\label{alg:iteration-gifs-d-quadtree}
    \State $T.n\leftarrow 0$ // initialize the number of triplets stored in $T$
    \State \textbf{create} root node of quadtree, $R$
    \For{$j\leftarrow 1$ \textbf{to} $L$}
        \For{$k\leftarrow 1$ \textbf{to} $D.n$}
            \For{$l\leftarrow 1$ \textbf{to} $D.n$}
                \State $(u,v)\leftarrow \phi_{j}(D[k].x,D[k].y,D[l].x,D[l].y)$
                \State \textbf{initialize} $r$ with some appropriate value
                \State \Call{Quadtree\_Search\_and\_Insert}{$R,u,v,r,T$}
            \EndFor
        \EndFor
    \EndFor
    \State $D \leftarrow T$
\EndFor
\EndFunction
\end{algorithmic}
\end{algorithm}

Once again, we will use two examples, defined below, to illustrate the functioning of the GIFS algorithms.
\begin{example}\label{example-filip-measures} %8.1
This example uses the GIFS $\mathcal{G}$ appearing in \cite[Example 16]{JMS16}. The IFS is defined by
\[
\left\{
  \begin{array}{ll}
       \phi_1((x_1,y_1),(x_2,y_2))=&(0.25x_1+0.2y_2,0.25y_1+0.2y_2) \\
       \phi_2((x_1,y_1),(x_2,y_2))=&(0.25x_1+0.2x_2,0.25y_1+0.1y_2+0.5)\\
       \phi_3((x_1,y_1),(x_2,y_2))=&(0.25x_1+0.1x_2+0.5,0.25y_1+0.2y_2))
  \end{array}
\right.
\]
on the region $X=[0,1]^2$, with $q_{1}=-2$, $q_{2}=0$ and $q_{3}=0$.
\end{example}

\begin{example}\label{example-final} %8.2
This example comes from \cite[Example 11.6]{COS21}. The IFS is defined by
\[
\left\{
  \begin{array}{ll}
    \phi_1((x_1,y_1),(x_2,y_2))=&(0.2x_1+0.25x_2+0.04y_2,0.16y_1-0.14x_2+0.20y_2+1.3) \\
    \phi_2((x_1,y_1),(x_2,y_2))=&(0.2x_1-0.15y_1-0.21x_2+0.15y_2+1.3,\\&0.25x_1+0.15y_1+0.25x_2+0.17) \\
    \phi_3((x_1,y_1),(x_2,y_2))=&(0.355x_1+0.355y_ 1+0.378,\\&-0.355x_1+0.355y_1+0.434-0.03y_2)
  \end{array}
\right.
\]
on the region $X=[-0.1,2.1]^2$, with $q_{1}=-1$, $q_{2}=0$ and $q_{3}=-7$.
\end{example}

In Table \ref{tab:gifs-d-quadtree} we present the execution times obtained with our \textsc{Fortran 2003} implementations of Algorithm \ref{alg:gifs-d} and Algorithm \ref{alg:gifs-d-quadtree}. The number of points generated after $4$ iterations was $n=2011229$ for Example \ref{example-filip-measures} and $n=13994321$ for Example \ref{example-final}.

Note that this enormous amount of points in the latter made us being unable to run Algorithm \ref{alg:gifs-d} in less than 12 hours of execution time, whereas with Algorithm \ref{alg:gifs-d-quadtree} it is possible to quickly obtain a solution. Again, we notice that $n_{\max}=64$ provides the smallest execution time for both examples.
\begin{table}
\centering
\caption{Maximum quadtree height and execution times of Algorithm \ref{alg:gifs-d-quadtree}.}\label{tab:gifs-d-quadtree}
\begin{tabular}{c|c|rrr}
  \hline
          & Algorithm \ref{alg:gifs-d} & \multicolumn{3}{c}{Algorithm \ref{alg:gifs-d-quadtree}} \\
     Example & Time [s] & $n_{\max}$ & $h$ & Time [s] \\\hline
  \ref{example-filip-measures} & 3055.7301 &    2 & 27 & 5.3230 \\
             &         &    4 & 27 & 3.8130 \\
             &         &    8 & 17 & 3.3790 \\
             &         &   16 & 14 & 3.1530 \\
             &         &   32 & 12 & 3.0860 \\
             &         &   64 & 11 & 3.0600 \\
             &         &  128 & 11 & 3.4380 \\
             &         &  256 & 10 & 3.5110 \\
             &         &  512 &  9 & 4.1260 \\
             &         & 1024 &  1 & 5.4640 \\
  \ref{example-final} & N/A &    2 & 26 & 155.6610 \\
                &         &    4 & 25 & 19.0220 \\
                &         &    8 & 24 & 16.1720 \\
                &         &   16 & 23 & 13.6780 \\
                &         &   32 & 22 & 13.1010 \\
                &         &   64 & 21 & 12.7690 \\
                &         &  128 & 20 & 13.8220 \\
                &         &  256 & 19 & 14.3010 \\
                &         &  512 & 17 & 16.4830 \\
                &         & 1024 & 13 & 21.0930 \\
  \hline
\end{tabular}
\end{table}

\begin{figure}
  \centering
  \begin{tabular}{c c}
  \includegraphics[width=3.85cm]{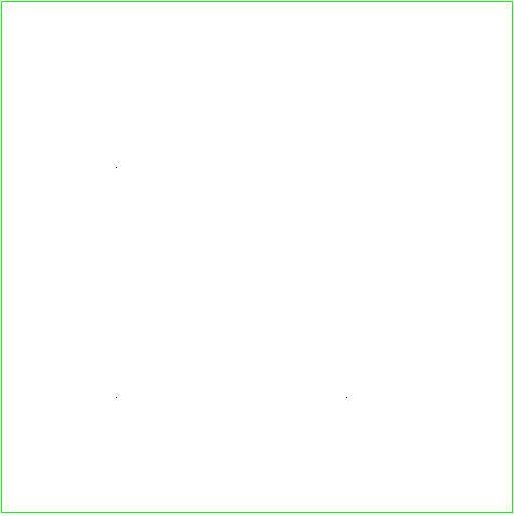}&\includegraphics[width=3.85cm]{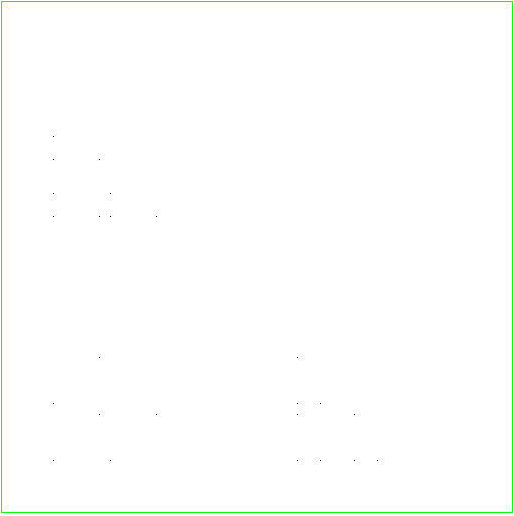}\\
  \includegraphics[width=3.85cm]{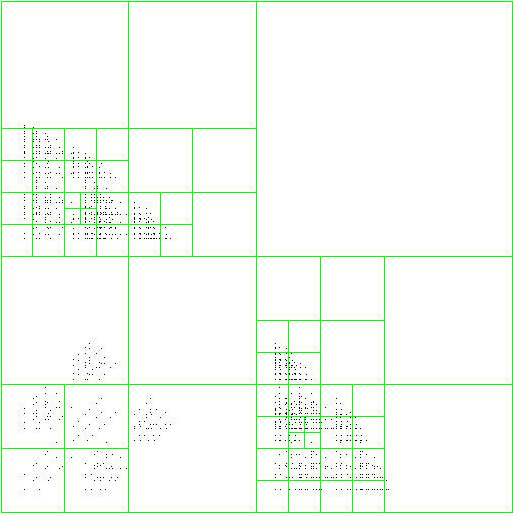}&\includegraphics[width=3.85cm]{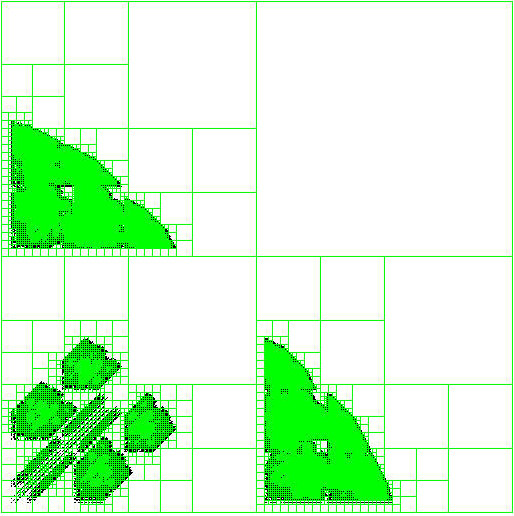}\\
  \end{tabular}
  \caption{The quadtrees placed over the fractal for Example \ref{example-filip-measures}, $n_{\max}=64$: from top to bottom, left to right, the images ($512\times 512$ pixels) show the points generated from iteration $1$ ($n=3$ points) to $4$ ($n=2011229$ points).}\label{fig:quadtree-exemplo-8-1}
\end{figure}

\begin{figure}
  \centering
  \begin{tabular}{c c}
  \includegraphics[width=3.85cm]{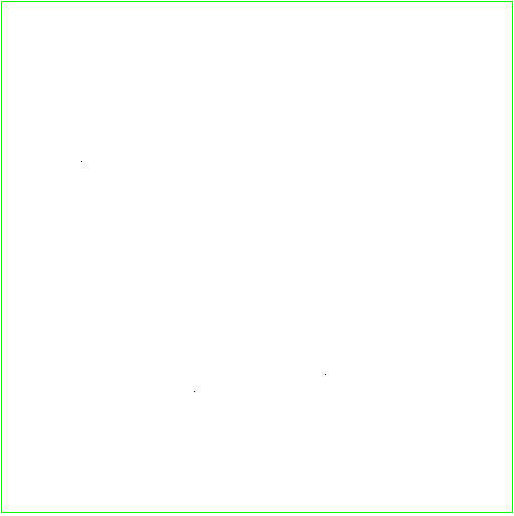}&\includegraphics[width=3.85cm]{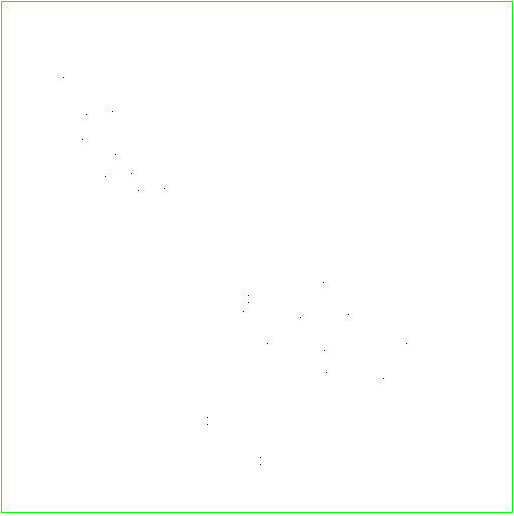}\\
  \includegraphics[width=3.85cm]{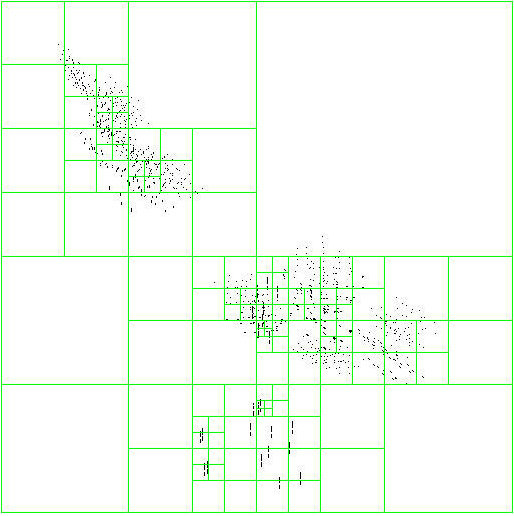}&\includegraphics[width=3.85cm]{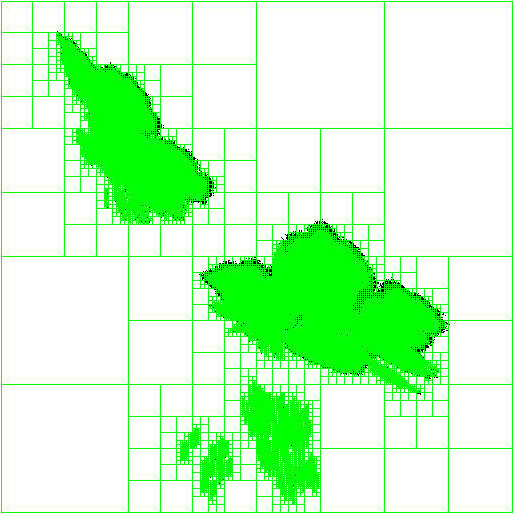}\\
  \end{tabular}
  \caption{The quadtrees placed over the fractal for Example \ref{example-final}, $n_{\max}=64$: from top to bottom, left to right, the images ($512\times 512$ pixels) show the points generated from iteration $1$ ($n=3$ points) to $4$ ($n=13994321$ points).}\label{fig:quadtree-exemplo-8-2}
\end{figure}

\begin{figure}
  \centering
  \begin{tabular}{c c}
  \includegraphics[height=3.85cm]{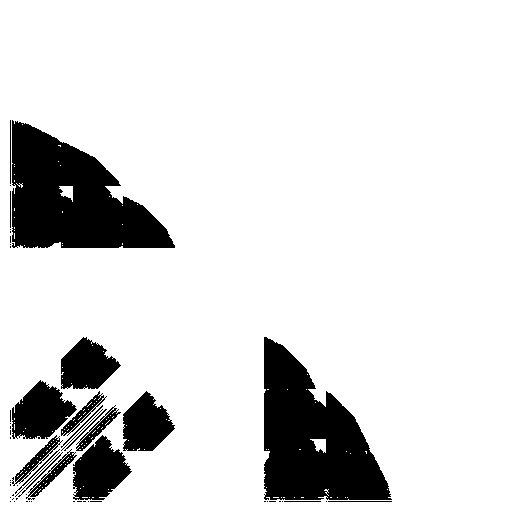}&\includegraphics[height=3.85cm]{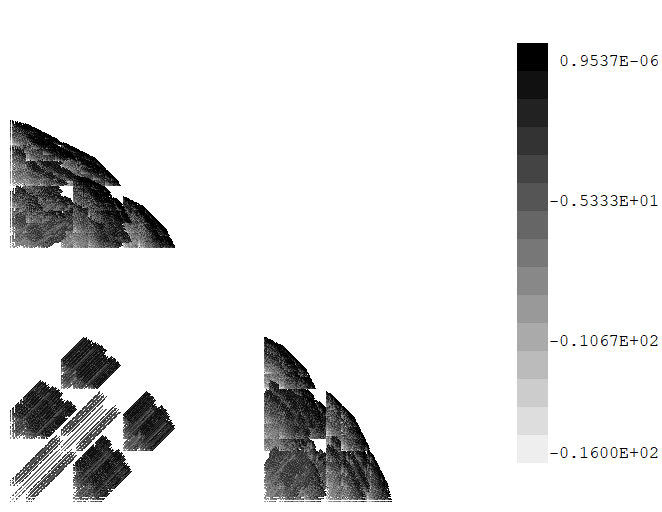}\\
  \includegraphics[height=3.85cm]{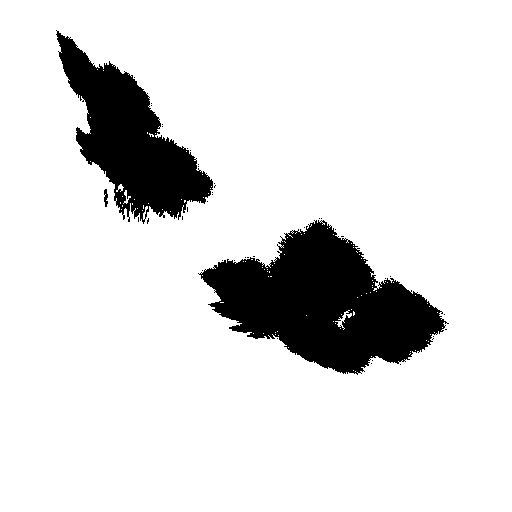}&\includegraphics[height=3.85cm]{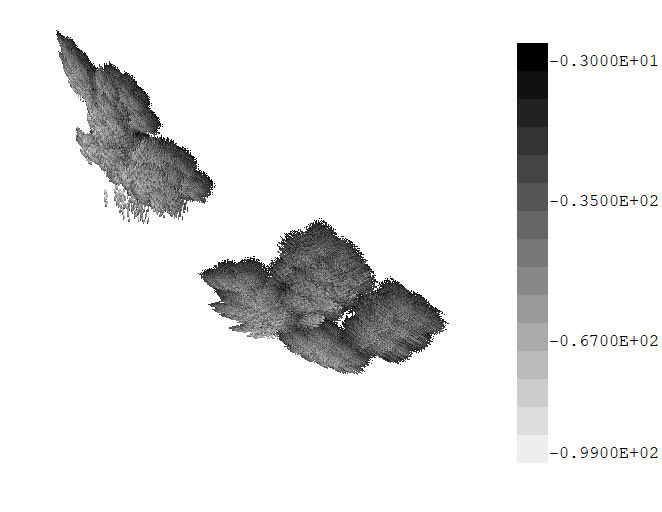}\\
  \end{tabular}
  \caption{The approximation of the attractor and the greyscale representation of the invariant measure for Example \ref{example-filip-measures} (top) and Example \ref{example-final} (bottom) idempotent GIFSs.}\label{fig:gifs-idempotent}
\end{figure}

As in Section \ref{sec:ifs-d-quadtree-in-practice}, we present in Figure \ref{fig:gifs-idempotent} the approximation of the attractor and the greyscale image representing the invariant measure for idempotent GIFSs.

\subsection{Complexity analysis}
\label{sec:complexity-gifs-d}
Under the same hypotheses assumed for Algorithm \ref{alg:ifs-d}, the total number of searches $S_{GIFS}$ made in Algorithm \ref{alg:gifs-d} is given by
\begin{equation}
S_{GIFS}=\sum_{i=1}^{N}{\sum_{j=1}^{L}{\sum_{k=1}^{n_{i-1}}{\sum_{l=1}^{n_{i-1}}{f_{k}}}}}
\label{eq:number-searches-gifs-d}
\end{equation}
and, assuming linear searches are used, it reduces to
\begin{equation}
S_{GIFS}=\frac{L^{2}-L^{2}L^{N}}{(L-1)^{2}}+\frac{1}{L-1}\sum_{i=1}^{N}{L^{i+1}\,\left(L^{L^{i-1}}\right)}
\end{equation}
and, since $L^{N+1}\,\left(L^{L^{N-1}}\right)$ dominates the summation, we may write
\begin{equation}
S_{GIFS}\simeq\frac{L^{2}-L^{2}L^{N}}{(L-1)^{2}}+
\frac{1}{L-1}\left(L^{N+1}\,L^{L^{N-1}}\right)\in O\left(L^{L^{N}+N}\right).
\label{eq:gifs-linear-search}
\end{equation}

For Algorithm \ref{alg:gifs-d-quadtree}, the total number of searches made, assuming $n_{\max}=64$ as before (see Section \ref{sec:complexity-ifs-d}), is given by
\begin{equation}
S_{GIFSq}=\frac{64}{L^{2}-1}\left(\frac{\left(L^{2}\right)^{N+1}-L^{2}}{L}\right)\in O(L^{2N})
\label{eq:B-gifs-quadtree}
\end{equation}
and, therefore,
\begin{equation}
S_{GIFSq}\ll S_{GIFS}.
\label{eq:comparison-gifs-gifsq}
\end{equation}

\section{Concluding remarks}
\label{sec:conclusion}
We have presented a description of the deterministic algorithm used to compute approximations of invariant measures and its attractors for IFS and GIFS, as well as a quadtree-based search algorithm that allows the use of these (G)IFS in a reasonable running time. The results presented show that our approach is effective in turning the deterministic algorithms for G(IFS) tractable.

\section*{Funding}
This research received no specific grant from any funding agency in the public, commercial, or not-for-profit sectors.

\section*{Conflict of interest}
The authors declare that they have no conflict of interest.


\end{document}